\documentclass[letterpaper, 10 pt, conference]{ieeeconf}  

\IEEEoverridecommandlockouts                              

\usepackage{amsmath} 

\usepackage{pgfplots}
\usepackage{tikz}
\usetikzlibrary{arrows.meta} 
\usepackage{amssymb}
\DeclareMathAlphabet{\mathcal}{OMS}{cmsy}{m}{n}
\usepackage{import}

\usepackage{multirow}
\usepackage{siunitx} 
\usepackage{booktabs}
\usepackage{tikz}
\usepackage{pgfplots}

\title{\LARGE \bf
	Data-driven control and transfer learning using neural canonical control structures*
}

\author{Lukas Ecker$^{1}$ and Markus Sch\"oberl$^{1}$
	\thanks{*This work has been supported by the COMET-K2 Center of the	Linz Center of Mechatronics (LCM) funded by the Austrian federal	government and the federal state of Upper Austria.}
	\thanks{$^{1}$Lukas Ecker and Markus Sch\"oberl are with Institute of Automatic Control and Control Systems Technology, Johannes Kepler University Linz, Altenberger Str. 69, A-4040 Linz, Austria (e-mail: {\tt\small \{lukas.ecker,markus.schoeberl\}@jku.at})}%
}

\begin{document}

\maketitle
\thispagestyle{empty}
\pagestyle{empty}

\begin{abstract}

An indirect data-driven control and transfer learning approach based on a data-driven feedback linearization with neural canonical control structures is proposed. 
An artificial neural network auto-encoder structure trained on recorded sensor data is used to approximate state and input transformations for the identification of the sampled-data system in Brunovsky canonical form. 
The identified transformations, together with a designed trajectory controller, can be transferred to a system with varied parameters, where the neural network weights are adapted using newly collected recordings. 
The proposed approach is demonstrated using an academic and an industrially motivated example.


%
\end{abstract}

\section{Introduction}

In recent years, data-driven control and transfer learning have emerged as promising techniques for controlling industrial systems, see \cite{9328227}.
Data-driven control distinguishes between direct and indirect approaches, see \cite{Krishnan2021}, where in direct data-driven control the control algorithm is derived directly from the data, while in indirect data-driven control a model is determined from the data and the control algorithm is designed based on the determined model.
Transfer learning is a technique that allows a model or controller trained on one task to be adapted and used for another related task, as shown in a robotic setting in \cite{Barrett2010TransferLF}.
These data-driven methods have the basic goal of reducing the modeling effort associated with control approaches that rely on a model.  One very popular nonlinear control method that heavily relies on an accurate model is  feedback linearization, c.f., \cite{isidori1995nonlinear}.
However, feedback linearization has also gained attention in a fully data-based scenario, see, e.g., \cite{YESILDIREK19951659},  
where the linearizing control law was functional approximated by two artificial neural networks \cite{Suykens1995ArtificialNN}.
Recent related works include, for example, \cite{Pasqualetti}, where a data-driven methodology for feedback linearization in a Koopman Operator setting is presented, \cite{westenbroek2020feedback}, where a reinforcement learning feedback linearization approach has been proposed, and \cite{Zahin}, where a linearization with nonlinear auto-regressive moving average models has been introduced.
In \cite{BruntonAutoEncoder}, an identification approach with a modified neural auto-encoder structure was used to embed nonlinear dynamics into a higher-dimensional linear Koopman framework. A data-based approach for feedback linearizable systems using Willems' fundamental lemma has been proposed in \cite{Alsalti2021DataBasedSA}.



This paper proposes an indirect data-driven feedback linearization approach, where the necessary transformations in canonical control form are implemented as an artificial neural network auto-encoder structure.
Approximating the discrete-time system in canonical Brunovsky form, where the sampling time is defined by the recorded data set, eliminates the need for a quasi-continuous implementation. The functional approximation of the transformation leads to an indirect approach that allows offline verification of the model in combination with a controller designed for the Brunovsky system. 
Furthermore, a validation of the model based on newly acquired sensor data is possible at any time, which is interesting for time-varying systems.
The structure of the discrete-time Brunovsky form, which is characterized by an inherent shift register chain, allows for straightforward planning and stabilization of desired trajectories.
When the controller is transferred to a non-nominal system, newly collected recordings are used to adjust the neural network weights of the feedback linearization until a Brunovsky behavior is re-established.
%
%

{
Besides a nonlinear academic example, the proposed data-driven control strategy is validated in simulation on a pair of elastic single mast stacker cranes, which are employed in industrial environments for material handling and storage. Note that several model-based controller and observer strategies have been already developed, see, e.g., \cite{EckerMalzerWahrburg} and \cite{7954034}.
%
%
Contrary to the above references, here the linear model of the single mast stacker crane without lifting unit is considered.
Nevertheless, the controller design for a series of cranes that can handle mechanical variations and disturbances while maintaining a high level of control performance without the need for manual adjustments can be challenging.
The proposed strategy aims to overcome this challenge by training a control model on a nominal single-mast stacker crane and then transferring the learned control to a deviating device. Newly recorded sensor data are then used to modify the controller until similar nominal behavior is achieved.}
%

The paper is organized as follows: In the problem statement in Section 2, the ideas of the control strategy and the data-driven framework are outlined. Section 3 addresses the assignment from a control theoretical point of view  whereas Section 4 outlines the data-driven control and transfer learning setting as well as the implementation and training of the neural canonical control structures. 
Section 5 presents the results of the data-driven controller design method using an academic case study and an industrially motivated application on single-mast stacker cranes.

\section{Problem Statement}


The control problem and its realization in a data-driven framework together with transfer learning are discussed.

\subsection{Control Strategy}
The proposed controller design aims to represent the system under consideration in a canonical control form and apply linear control theory methods for disturbance suppression and trajectory stabilization. The controller structure is composed of an inner-loop and an outer-loop control. 
The inner loop, consisting of a feedback linearization, converts the model into an equivalent linear control system in Brunovsky canonical form, see \cite{isidori1995nonlinear}, using state and input transformations.
The outer control-loop, which relies on the resulting linear canonical control system of the inner-loop, is used for tracking control and disturbance suppression.
While the outer controller is designed based on the Brunovsky structure, the feedback linearization in the inner loop is to be derived from historical sensor data rather than from a mathematical model. 
The data-based approach with sampled data naturally leads to a discrete-time system  feedback linearization as discussed in \cite{Grizzle} and  \cite{Nijmeijer_Nonlinear_dynamical_control_systems}.
The control strategy combined with a data-driven approximation of the linearizing transformations has several advantages, including trivial computation of target trajectories, stability and performance in the presence of disturbances and parameter variations, and the ability to transfer similar behavior of a nominal system to a deviating system by adopting the outer linear controller.



\subsection{Data-Driven Control and Transfer Learning}

The main objective is to reduce the mathematical modeling and controller design effort while allowing the adaptation and transfer of the control model to other non-nominal systems. Along with this, similar control characteristics and behaviors should be maintained.
The proposed indirect inner-loop data-driven feedback-linearization approximates the transformations to Brunovsky canonical form with historical sensor data.
In contrast to e.g. \cite{YESILDIREK19951659} where a linearizing control law is immediately approximated, here the state and input transformations are learned separately, resulting in a data-driven identification of a sampled-data system. 
Thus, the indirect approach leads to an approximation of a continuous-time physical process $\dot{x} = f_c(x,u)$
as discrete time-system 
\begin{align}\label{equ:sampled_data}
	x_{k+1} = f_d(x_k,u_k),
\end{align}
where $k$ denotes the time index, i.e., $x(kT_a) = x_k$, with the sampling time $T_a$  defined by the recorded sensor data.
The successor state defined by the state $x_k \in \mathcal{D}_x \subset \mathbb{R}^{n}$, the input $u_k \in \mathcal{D}_u \subset \mathbb{R}$, and the transition mapping (\ref{equ:sampled_data}) is abbreviated as $x_k^+ = x_{k+1} = f_d(x_k,u_k)$.
The data-driven modeling approach uses $N_s$ samples of the state $\mathcal{X} = \{x_i\}_{i=0}^{N_s}$, the input $\mathcal{U} = \{u_i\}_{i=0}^{N_s}$, and the successor state $\mathcal{X}^+ = \{x^+_i\}_{i=0}^{N_s}$ from previously recorded experiments.
%
%
The superimposed outer-loop trajectory controller is not derived from data, but is designed on the Brunovsky system using methods from linear control theory, such as pole placement.
When the control model is transferred to a non-nominal system, the transformations must be readjusted to maintain the Brunovsky form in the inner loop. Data sets from the transferred system with the nominal controller are used to adapt the nominal network to the modified system.

\section{Control Structure and Strategy}

This section recapitulates ideas from nonlinear control theory for transforming the system into the Brunovsky canonical form.  It provides a summary of the methods used to generate desired trajectories and to design the outer-loop controller.


\subsection{Transformation in Brunovsky Canonical Form}

\begin{figure}
	\begin{center}
		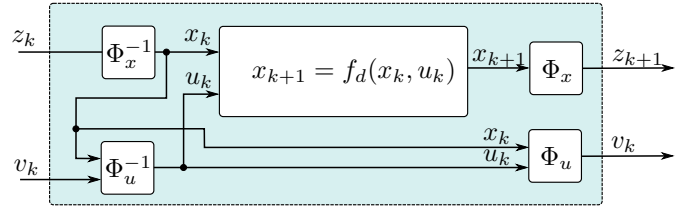
		\caption{Schematic conversion of the sampled data system into the canonical Brunovsky form according to (\ref{equ:sampled_data_system_to_brunovsky}).}
		\label{fig:canonical_structures_box_sole_inverse}
	\end{center}
\end{figure}

Considered are systems that can be transformed with a proper state transformation $z_k = \Phi_x(x_k) : \mathcal{D}_x \rightarrow \Phi_x(\mathcal{D}_x)$ and an input transformation $v_k = \Phi_u(x_k,u_k): \mathcal{D}_x \times \mathcal{D}_u \rightarrow \Phi_u(\mathcal{D}_x,\mathcal{D}_u)$ according to
\begin{align}
	z_{k+1} &= \Phi_x \circ f_d(\Phi_x^{-1} (z_k),\Phi_u^{-1}(\Phi_x^{-1}(z_k),v_k))
	\label{equ:sampled_data_system_to_brunovsky}
\end{align}
into a linear discrete-time system with new state $z_k \in \Phi_x(\mathcal{D}_x)$ in Brunovsky canonical form
\begin{align}
	\label{equ:brunovsky_system}
z_{k+1} = \begin{bmatrix}
	0 & 1 & \cdots & 0 \\
	\vdots & \vdots & \ddots & \vdots \\
	0 & 0 & \cdots & 1 \\
	0 & 0 & 0 & 0 \\
\end{bmatrix} z_k + 
\begin{bmatrix}
	0\\
	\vdots\\
	0\\
	1\\
\end{bmatrix}
v_k.
\end{align}
This conversion is illustrated schematically in Fig. \ref{fig:canonical_structures_box_sole_inverse}.
Note that the proposed approach in this work focuses on the input state linearizable scenario with $\dim(z_k) = \dim(x_k) = n$ and single input $\dim(u_k) = 1$.
In contrast to the continuous-time scenario, where the Brunovsky form corresponds to an integrator chain, the behavior in the discrete-time scenario is equivalent to a shift register, and thus abbreviated as $ z_{k+1} = A_B z_k + b_B v_k = \sigma(z_k,v_k)$. Due to the bijective assumptions on the transformations, the system in original coordinates $x_{k+1} = f_d(x_k,u_k)$ can be expressed by means of the canonical shift system and the transformations as 
\begin{equation}
	x_{k+1} = \Phi_x^{-1} \circ \sigma(\Phi_x(x_k), \Phi_u(x_k,u_k)). \label{equ:brunvosky_to_sampled_data}
\end{equation}

\begin{figure}
	\begin{center}
		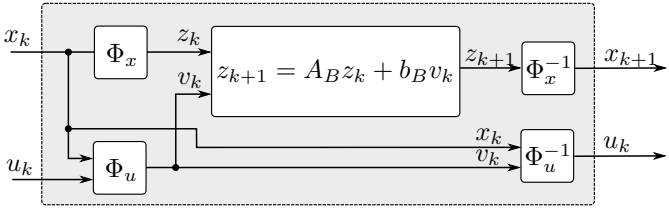
		\caption{Representation of the sampled data system as transformed linear system in Brunovsky canonical form according to (\ref{equ:brunvosky_to_sampled_data}).} 
		\label{fig:representation_dynamical_system_as_canonical_control_form}
	\end{center}
\end{figure}

A schematically interpretation of this equation is given in Fig. \ref{fig:representation_dynamical_system_as_canonical_control_form}. 
Since the Brunovsky structure is uniquely defined by the dimension of the state $z_k$, only the transformations $\Phi_x$ and $\Phi_u$ have to be determined for the approximation of the original system. 
The existence of these transformations can be traced back to a distribution test, see \cite{Nijmeijer_Nonlinear_dynamical_control_systems} Algorithm 14.2, where the subsequent determination of the transformation requires the solution of a set of ordinary differential equations.
The proposed approach seeks to overcome this burden by implementing the transformations as neural networks and to train them on historical sensor data.
This attempts to impose the behavior of the Brunovsky system on the original system.
%
The subsequent trajectory planning and the outer controller do not need the transformations in an analytical form, the implementation as neural networks is sufficient. 
\subsection{Trajectory Planning}\label{sec:trajectory_planning}
For a desired change of the operating point from an initial position $x_0$ to a final position $x_N$ in $N$ sampling steps, the trajectory is planned in Brunovsky coordinates.  For this purpose, the initial and target coordinates $z_0 = \Phi_x(x_0)$ and $z_N = \Phi_x(x_N)$ are determined using the state parametrization. The first Brunovsky state $y_1(k)$ is parameterized with an approach $y_{d,1}(k) = \lambda(k)^T \alpha$ linear in the coefficients $\alpha \in \mathbb{R}^{n_\alpha}$. The parameters $\lambda(k)^T = \left[\lambda_1(k), \ldots, \lambda_{n_\alpha}(k)\right]$ with $\lambda_i \in \mathbb{R}$ have to satisfy the initial and end conditions. A polynomial approach for lambda with $\lambda(k)^T = \left[1, k , \ldots, k^{n_\alpha}\right]$ and $n_\alpha = 2n-1$ was chosen. Due to the shift-register structure of the Brunovsky system, the remaining Brunovsky states are determined by $y_{d,i}(k) = \lambda(k+i-1)^T \alpha$ for $i=1,\ldots,n$ and thus
$$z_d(k) = \begin{bmatrix}
	y_{d,1}(k) \\
	\vdots\\
	y_{d,n}(k)
\end{bmatrix}
= \begin{bmatrix}
	\lambda(k)^T \\
	\vdots\\
	\lambda(k+n-1)^T
\end{bmatrix}\alpha = \mathcal{S}(k)\alpha.
$$
The unknown coefficients $\alpha$ can be determined by solving the linear equation system of the boundary conditions
$$
\begin{bmatrix}
	z_0\\
	z_N\\
\end{bmatrix} 
=
\begin{bmatrix}
	z_d(0)\\
	z_d(N)\\
\end{bmatrix}
=
\begin{bmatrix}
	\mathcal{S}(0)\\
	\mathcal{S}(N)\\
\end{bmatrix}
\alpha.
$$
The target input is defined by $v_d(k) = \lambda(k+n)^T\alpha$.

\subsection{Outer-Loop Controller}
\begin{figure}
	\begin{center}
		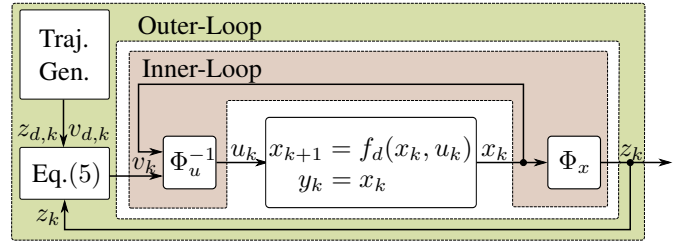
		\caption{Schematic of the control strategy with inner transformation to Brunovsky canonical form and an outer trajectory controller.} 
		\label{fig:controller_schematic}
	\end{center}
\end{figure}

The outer-loop controller is designed for the Brunovsky system, with the aim to impose an exponentially stable error dynamic 
\begin{equation*}
	\begin{bmatrix}
		{e}_{1,k+1}\\
		\vdots\\
		{e}_{n,k+1}
	\end{bmatrix} 
	=
	\begin{bmatrix}
		0 & 1 & \cdots & 0 \\
		\vdots & \vdots & \ddots & \vdots \\
		0 & 0 & \cdots & 1 \\
		-a_0 & -a_1 & \cdots & -a_{n-1} \\
	\end{bmatrix}
	\begin{bmatrix}
		{e}_{1,k}\\
		\vdots\\
		{e}_{n,k}
	\end{bmatrix} 
\end{equation*} 
for the trajectory error $e_1 = y_1 - y_{1,d}$ with the desired target trajectory from Section \ref{sec:trajectory_planning}. The corresponding control law can be deduced as
\begin{equation}
	v_k = - \left[ a_0, \ldots, a_{n-1}
	\right] \cdot 
\left(z_k - z_{d,k}\right)  + v_{d,k} 
\end{equation}
with proper chosen coefficients $a_j$, $j = 0,\ldots,n-1$. The coefficients correspond to the characteristic polynomial  of the error dynamics and can be determined, e.g., by pole placement. The overall control schematic divided in an inner- and outer-loop is presented in Fig. \ref{fig:controller_schematic}.

\section{Neural canonical control structures}
Instead of laboriously deriving the transformations $\Phi_x$ and $\Phi_u$ from a mathematical model, they are defined as artificial neural networks and trained using recorded sensor data.
The transformation should approximate the system schematically defined in Fig. 
\ref{fig:representation_dynamical_system_as_canonical_control_form}
as closely as possible according to the principal optimization problem
\begin{align*}
	\hspace{-0.5cm}
	\min_{\Phi_x,\Phi_x^{-1},\Phi_u,\Phi_u{-1}}  &\lVert x^+ - \Phi_x^{-1} \circ \sigma(\Phi_x(x),\Phi_u(x,u))\rVert  \\
	\text{s.t.} 
	\quad & x = \Phi_x^{-1}(\cdot) \circ \Phi_x(x) 
	\\
	\quad & u = \Phi_u^{-1}(x,\cdot) \circ \Phi_u(x,u) \\
	\quad &    x \in \mathcal{D}_x, \quad u \in \mathcal{D}_u.
\end{align*}
The invertibility of the transformations is crucial for the design, otherwise a representation of the system in canonical Brunovsky form and the feedback linearization will fail. This must also be taken into account when implementing the transformations as neural networks.
One solution would be to implement $\Phi_x$ as an analytically invertible neural network, but this restricts the approach to invertible activation functions, constant number of units in each layer $n_l$, and to the input-state linearizable scenario with $\dim(x_k) = \dim(z_k)$. 
Although this does not contradict the imposed assumptions, the planned extension of this data-driven control approach to flat systems in the discrete-time case, c.f.,  \cite{diwold2} and \cite{diwold1},  would not be possible. 
%
Therefore, $\Phi_x$ and $\Phi_u$ as well as $\Phi_x^{-1}$ and $\Phi_u^{-1}$ are implemented as separate artificial neural networks. 
As already indicated in Fig. \ref{fig:representation_dynamical_system_as_canonical_control_form}, this task is similar to learning an auto-encoder structure with $\Phi_x$, $\Phi_u$ as encoders and $\Phi_x^{-1}$, $\Phi_u^{-1}$ as decoders, where the encoded Brunovsky states are shifted through the Brunovsky system (\ref{equ:brunovsky_system}) and decoded back to the original coordinates.
In addition to achieving high prediction accuracy, the invertibility constraint for the two transformations $x = \Phi_x^{-1}(\cdot) \circ \Phi_x(x)$
and $u = \Phi_u^{-1}(x,\cdot) \circ \Phi_u(x,u)$ must be ensured.
%
As is common with auto-encoders, see \cite{BruntonAutoEncoder}, this was accomplished by including the constraints as penalty terms in the cost function along with the prediction loss.
%
%
The transformations are functional approximated as feed forward artificial neural networks. The numerical results in Section \ref{sec:experimental_results} are achieved with a one hidden layer network topology and nonlinear activation  function $\phi(s) = \left[ \phi_1(s_1),\ldots,\phi_{n_l}(s_{n_l})\right]^T \in \mathbb{R}^{n_l}, s_i \in \mathbb{R}$. 
The transformations were thus approximated as 
\begin{equation}
	\label{equ:transformation_nn}
\hat{\Phi}(x_k; W_{\hat{\Phi},i}, b_{\hat{\Phi},i}) = W_{\hat{\Phi},2} \phi(W_{\hat{\Phi},1} x_k + b_{\hat{\Phi},1}) + b_{\hat{\Phi},2}\end{equation} with the parameters to be optimized such as weights $W_{\hat{\Phi},1} \in \mathbb{R}^{n_l \times n}$, $W_{\hat{\Phi},2} \in \mathbb{R}^{n \times n_l}$ and biases $b_{\hat{\Phi},1} \in \mathbb{R}^{n_l}, b_{\hat{\Phi},2}\in \mathbb{R}^{n}$.
%
The complete loss function for training the auto-encoder structure can be synthesized as
\begin{align*}
	L =  \alpha_1 L_{\text{rec},x} + \alpha_2 L_{\text{rec},u} + \alpha_3 (L_{\text{pred},1}+  L_{\text{pred},2})
\end{align*}
with the state reconstruction loss  (\ref{equ:rec_x}), the input reconstruction loss (\ref{equ:rec_u}), the prediction losses (\ref{equ:pred1}) and (\ref{equ:pred2}), and the hyperparameters $\alpha_1,\alpha_2,\alpha_3 \in \mathbb{R}$.
\begin{align}
	L_{\text{rec},x} &= || x - \Phi_x^{-1} \circ \Phi_x (x) ||_{\text{MSE}} \label{equ:rec_x}\\
	L_{\text{rec},u} &= || u - \Phi_u^{-1}\circ(x, \Phi_u (x, u)) ||_{\text{MSE}} \label{equ:rec_u}\\
	L_{\text{pred},1} &= || x^+ - \Phi_x^{-1} \circ \sigma(\Phi(x),\Phi_u(x,u))||_{\text{MSE}} \label{equ:pred1}\\
	L_{\text{pred},2} &= || \Phi_x(x^+) - \sigma(\Phi(x),\Phi_u(x,u))||_\text{MSE}\label{equ:pred2}
\end{align}


\section{Experimental Results}\label{sec:experimental_results}

The proposed indirect data-driven approach is demonstrated by two examples. The first  example (\ref{sec:academic_example}) considers an academic 3-dimensional nonlinear system. Prediction and reconstruction performances are compared to illustrate the principle capabilities of the approach.
The second example (\ref{sec:industrial_example}) shows the approach in an industrially motivated application with combination of trajectory control and transfer learning.
The auto-encoder structures were trained in the open source machine learning framework PyTorch with the Adam optimization algorithm.

\subsection{Academic Example}
\label{sec:academic_example}

The proposed data-driven approach is applied to the 3-dimensional nonlinear academic example
\begin{equation}
	x_{k+1} = \begin{bmatrix}
		\frac{x_{2,k}}{x_{1,k}+2} \\
		\frac{x_{2,k} x_{3,k}}{x_{1,k}+2} + 2 x_{3,k}\\
		\frac{u_{1,k}}{x_{2,k}+2}
	\end{bmatrix},
\end{equation}  which is known to be exactly linearizable by static feedback according to the theory in \cite{Grizzle}.
The feedback-linearization transformations were approximated according to (\ref{equ:transformation_nn}) with the sigmoid activation function  $\phi_i(s_i) = \frac{1}{1+e^{-s_i}},s_i \in\mathbb{R}$ and $n_l = 80$ units in the hidden layer.
The neural networks were trained with $N_S = 50000$ recorded samples of randomly initialized trajectories and random input signals. The prediction and reconstruction results of a validation set after 10000 epochs of training are presented in Fig. \ref{fig:academic_example}. The left column shows the very high accuracy of the identified system with internal representation in Brunovsky form. The right column shows the high reconstruction accuracy of the auto-encoder structure required for the proposed control strategy. 
This example is intended to demonstrate the basic capabilities of the algorithm and not to make a qualitative statement. No tuning of the hyperparameters was performed and the comparison of different numbers of neurons was omitted due to space constraints.
{ Nevertheless, very good results could be obtained in this academic example, which should demonstrate the basic performance of the approach.}

\DeclareRobustCommand\full  {\tikz[baseline=-0.6ex]\draw[thick,orange,line width = 0.5mm] (0,0)--(0.5,0);}
\DeclareRobustCommand\dashed{\tikz[baseline=-0.6ex]\draw[thick,dashed, black,style=dashed,line width = 0.5mm] (0,0)--(0.54,0);}

\DeclareRobustCommand\fullblack  {\tikz[baseline=-0.6ex]\draw[thick,black,line width = 0.5mm] (0,0)--(0.5,0);}

\begin{figure}[t]
	\begin{center}
		\hspace{-0.5cm}
		\includegraphics{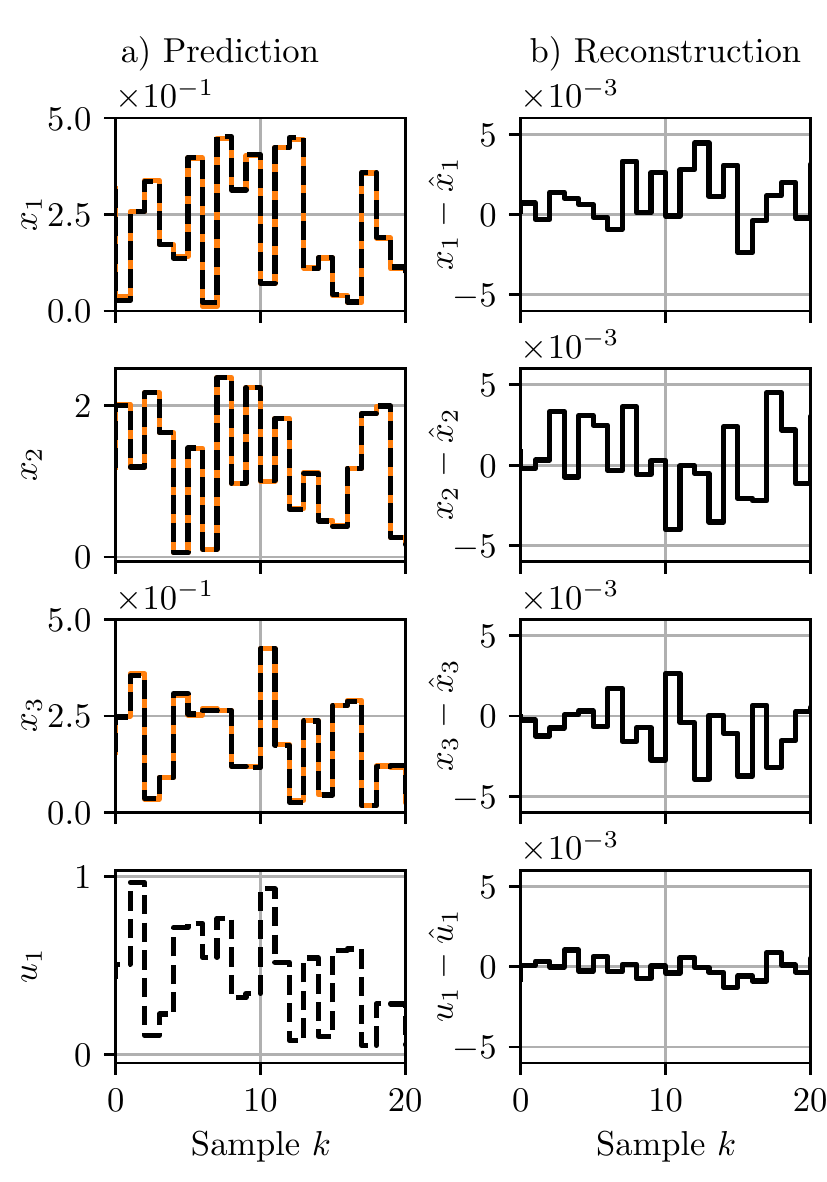}
		\caption{a) Prediction: Comparison between validation trajectory (\dashed) and identified model (\full) with feed forward input signal and same initial conditions.
		b) Reconstruction: Error (\fullblack) between validation trajectory $x_k$, $u_k$ to its  auto-encoder reconstruction $\hat{x}_k = \Phi_x^{-1}( \Phi_x(x_k))$, $\hat{u}_k = \Phi_u^{-1}(x_k, \Phi_u(x_k,u_k))$.}
		\label{fig:academic_example}
	\end{center}
\end{figure}

\subsection{Industrially Motivated Application}
\label{sec:industrial_example}

The data-driven design process is now demonstrated for a single mast stacker crane without lifting unit. A brief description of the model is given in Section \ref{sec:mathematical_model}.
The design process is divided into two parts.
The first part, as outlined in Section \ref{sec:nominal_controller_design}, presents the results of the data-driven controller for a nominal crane $\Sigma_N$ that is designed to suppress disturbances and to position the crane along planned and stabilized trajectories. 
The second part, covered in Section \ref{sec:transfer_controller_design}, illustrates the application of the nominal controller to a deviating crane $\Sigma_T$ as well as the transfer learning and adaptation of the controller.
%
%

\subsubsection{Mathematical Model}
\label{sec:mathematical_model}
%
%
%

\begin{table}[t]
	\begin{center}
			
			\caption{Parameters of the stacker crane models.}\label{tab:parameter}
			{
					\begin{tabular}{c|ccccc}
							\toprule
							
							System & $L$ / \SI{}{\meter} & $m_c$ / \SI{}{\kilogram} & $m_h$ / \SI{}{\kilogram} & $\rho A$ / \SI{}{\kilogram/\meter} & $EI$ / \SI{}{\newton \square \meter}\\
							\midrule
							\midrule 
							$\Sigma_N$ & 0.53 &  13.10 & 0.32 & 2.10  & 14.97 \\
							\midrule 
							$\Sigma_T$ & 0.53 &  12.72 & 0.34 & 2.26  & 14.28 \\
							\bottomrule
					\end{tabular}}
		\end{center}
	
\end{table}

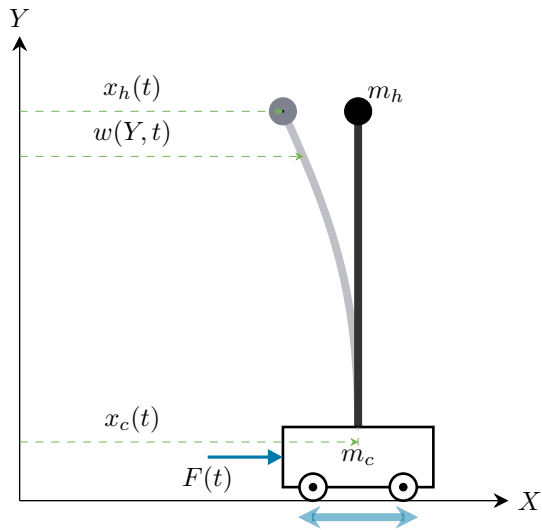
\begin{figure}[t]
	\centering
	\definecolor{jkuBlue}{RGB}{4,110,152}
\definecolor{jkuBlue}{RGB}{0,120,170}
\definecolor{jkuCyan}{RGB}{100,180,190}
\definecolor{jkuYellow}{RGB}{230,195,35}
\definecolor{jkuGrey}{RGB}{125,130,140}
\definecolor{jkuDarkGrey}{RGB}{51,51,51}
\definecolor{jkuLightGreen}{RGB}{195,215,75}
\definecolor{jkuGreen}{RGB}{115,180,85}
\definecolor{jkuPurple}{RGB}{145,75,130}
\definecolor{jkuRed}{RGB}{205,90,80}

\begin{tikzpicture}

\draw[line width=3,jkuGrey,opacity=0.5] (1,0.8) to [out=90, in=-65] (0,5);

\draw[line width=3,jkuDarkGrey] (1,0.8) rectangle (1,5);

\filldraw[fill=black] (1,5) circle (0.18);
\node [above right] at (1,5) {$m_h$};
\filldraw[fill=jkuDarkGrey,jkuGrey] (0,5) circle (0.18);
\filldraw[fill=jkuGrey] (0,5) circle (0.01);


\draw[line width=1] (0,0) rectangle (2,0.8);
\node at (1,0.4) {$m_c$};
\filldraw[fill=white,line width=1] (0.4,0) circle (0.18);
\filldraw[fill=white,line width=1] (1.6,0) circle (0.18);
\filldraw[fill=black] (0.4,0) circle (0.05);
\filldraw[fill=black] (1.6,0) circle (0.05);

\draw[-{Triangle[length=2mm,width=2mm]},line width=1.2,jkuBlue] (-1,0.4) -- (0,0.4);
\node [left,below] at (-1,0.4) {$F(t)$};

\draw[-{Stealth[length=2mm,width=2mm]},line width=0.5] (-3.5,-0.18) -- (3.0,-0.18);
\node [right]  at (3.0,-0.18) {$X$};

\draw[-{Stealth[length=2mm,width=2mm]},line width=0.5] (-3.5,-0.18) -- (-3.5,6);
\node [above]  at (-3.5,6) {$Y$};

\draw[-{Stealth[length=1mm,width=1mm]},line width=0.25,jkuGreen,dashed] (-3.5,5) --   (0,5); 
\draw[-{Stealth[length=1mm,width=1mm]},line width=0.25,jkuGreen,dashed] (-3.5,4.4) -- (0.27,4.4);
\node [above]  at (-2,5.0) {$x_h(t)$};
\node [above]  at (-2,4.4) {$w(Y,t)$};


\draw[-{Stealth[length=1mm,width=1mm]},line width=0.25,jkuGreen,dashed] (-3.5,0.6) -- (1,0.6);
\node [above]  at (-2,0.6) {$x_c(t)$};
\draw[line width=0.25,jkuGreen,dashed] (1,0.55) -- (1,0.8);

\draw[{Stealth[length=3mm,width=3mm]}-{Stealth[length=3mm,width=3mm]},line width=3,jkuBlue,opacity=0.5] (1-0.8,-0.4) -- (1+0.8,-0.4);

\end{tikzpicture}
	\caption{Schematic of the Single Mast Stacker Crane.}
	\label{fig:Testbench_Model}	
\end{figure}

The principal schematic of the considered single mast stacker crane without lifting unit is depicted in Fig. \ref{fig:Testbench_Model}.
The position of the rigid driving cart with mass $m_c$ is denoted by $x_c$ and horizontal position of the tip mass $m_h$ is specified by the coordinate $x_h$.
When the lifting unit is neglected, the problem can be described by a linear system of ordinary and partial differential equations.
A detailed derivation of the mathematical model using a variational approach can be found in \cite{EckerMalzerWahrburg} and \cite{7954034}. The governing equations include the partial differential equation of the Euler-Bernoulli beam, the momentum equations of the driving unit and tip mass, as well as the boundary constraints.
By approximating the beam deformation $w(Y,t)$ with a first-order Rayleigh-Ritz ansatz
$
w(Y,t) = x_c(t) + \Psi(Y) \bar{q}(t)
$
and proper spatial ansatz function $\Psi(Y)$, 
the Euler-Lagrange equations yield a finite-dimensional 
mechanical system
$
M \ddot{{q}} + C q = G u
$
with generalized coordinates $ q = [x_c, \bar{q}]^T$ and input force $u = F$. The matrices are
$$
M = 
\begin{bmatrix}
	m_{11} & m_{12} \\
	m_{12} & m_{22} 	
\end{bmatrix},\;\;
C= 
\begin{bmatrix}
	0 \\
	c_2	
\end{bmatrix}, \;\text{and}\;\;
G
=
\begin{bmatrix}
	1\\
	0
\end{bmatrix}
$$
with $m_{11} = \rho A L + m_c + m_h$, $m_{12} = m_h \Psi(L) + \rho A \int_{0}^{L}\Psi(Y)\text{d}Y$, $m_{22} = m_c \Psi(L)^2 + \rho A \int_{0}^{L} \Psi(Y)^2\text{d}Y$ and $c_2=  EI \int_{0}^{L} 
\left(\partial^2 \Psi(Y) / \partial Y^2 \right)^2
\text{d}Y $, where $L$ denotes the length of the elastic mast, $\rho A$ the mass density and $EI$ the bending stiffness.
This approximation of the single mast stacker crane serves as simulation benchmark example for the proposed data-driven controller design.



\subsubsection{Nominal Controller}\label{sec:nominal_controller_design}

%
The initial controller design for the nominal system $\Sigma_N$ is based on historical sensor data from approximately 400 recorded trajectories with overall $N_s=320000$ samples and a sampling time of $T_a = \SI{5}{\milli \second}$. 
These trajectories include stabilization of different operating points, changes of operating points and excitations with additional random inputs each using a trivial PD control law. Again, the sigmoid activation function was applied and $n_l = 120$ neurons in the hidden layer.
The results of the trained controller for the nominal system $\Sigma_N$ are presented in Fig. \ref{fig:smc_trajectory_01}. Shown is a trajectory planned according to the procedure in Section \ref{sec:trajectory_planning} for a working point change from $x_{c,0} = \SI{0}{\meter}$ to $x_{c,N} = \SI{1}{\meter}$ within $N = 400$ samples and an initial positioning error of \SI{0.1}{\meter}.
The initial trajectory error is corrected by the stabilizing outer loop controller with error dynamic poles tuned to a magnitude less than 1.
The system behavior can be easily adjusted by the error dynamics.
%
\DeclareRobustCommand\full  {\tikz[baseline=-0.6ex]\draw[thick,orange,line width = 0.5mm] (0,0)--(0.5,0);}
\DeclareRobustCommand\dashed{\tikz[baseline=-0.6ex]\draw[thick,dashed, black,line width = 0.5mm] (0,0)--(0.54,0);}

\begin{figure}[t]
	\begin{center}
		\hspace{-1.cm}
		\includegraphics{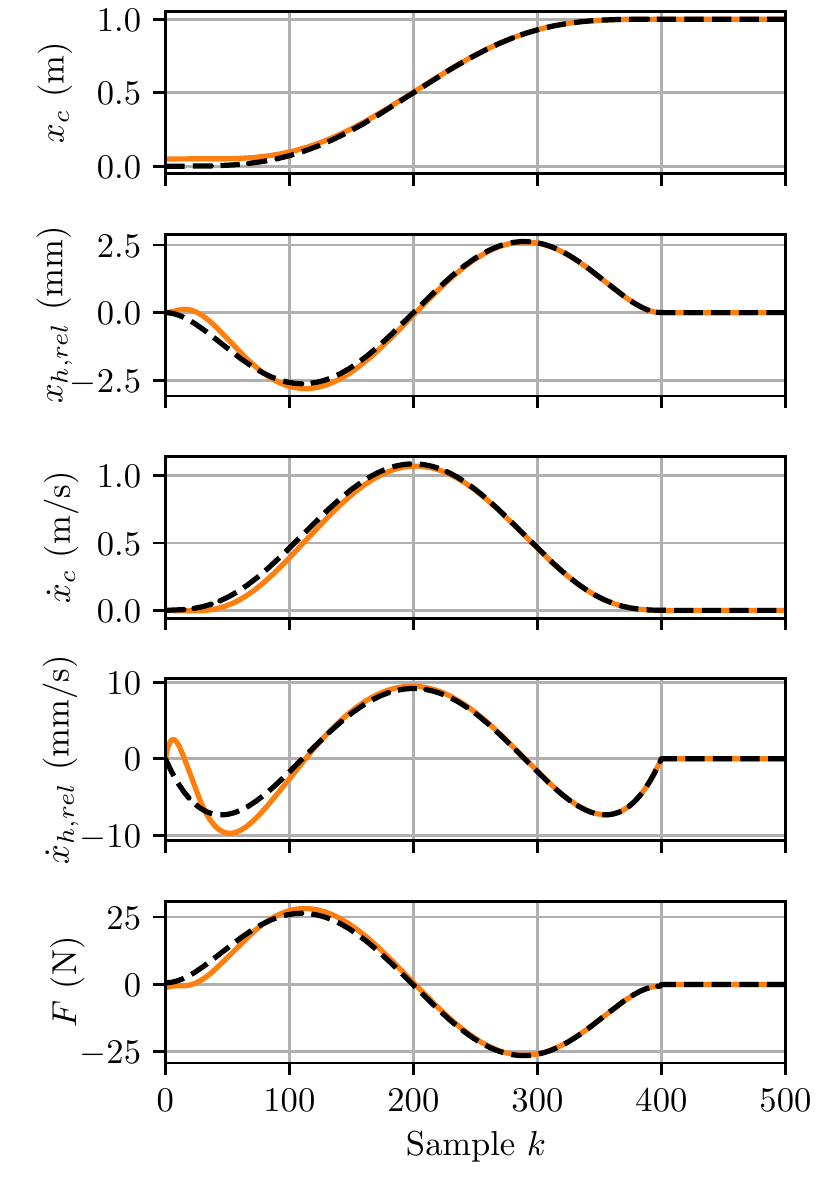}
		\caption{State and input trajectories of the controlled nominal system  $\Sigma_N \; (\full)$ in comparison to the desired trajectories (\dashed).} 
		\label{fig:smc_trajectory_01}
	\end{center}
\end{figure}

\subsubsection{Transferred Controller}\label{sec:transfer_controller_design}

The parameters of the target crane $\Sigma_T$ are listed in table \ref{tab:parameter}. The parameter deviations are within the range $\pm \SI{10}{\percent}$ of the nominal crane, except for the mast length. The nominal controller is applied to the target crane $\Sigma_T$ with the same objective of changing the operating position from $x_{c,0} = \SI{0}{\meter}$ to $x_{c,N} = \SI{1}{\meter}$. 
The state trajectories of the varied crane $\Sigma_T$ versus the nominal crane $\Sigma_N$ are shown in Fig. \ref{fig:smc_trajectory_04}.
The control error increased significantly compared to the nominal case and also an undesired oscillatory behavior appeared with the given parameter constellation.
Nevertheless, the controller manages the desired change in position. 
In order to increase the performance, sensor data of 50 experiments with the nominal controller were collected and used to adjust the neural network weights of the auto-encoder structure. 
%
%
The result after the transfer learning procedure is also shown in Fig. \ref{fig:smc_trajectory_04}.
The comparison between the nominal controller and the tuned controller on crane $\Sigma_T$ shows a significant improvement.
The adjusted controller behaves similarly to the nominal controller and also the introduced oscillations could be suppressed.
However, the small deviations from the nominal target trajectory are due to the adjustment of the auto-encoder structure, which of course also slightly changes the target trajectories.

\DeclareRobustCommand\full  {\tikz[baseline=-0.6ex]\draw[thick,orange,line width = 0.5mm] (0,0)--(0.5,0);}
\DeclareRobustCommand\dashed{\tikz[baseline=-0.6ex]\draw[thick,dashed, black,line width = 0.5mm] (0,0)--(0.54,0);}

\DeclareRobustCommand\fulla  {\tikz[baseline=-0.6ex]\draw[thick,orange,line width = 0.5mm] (0,0)--(0.5,0);}
\DeclareRobustCommand\fullb  {\tikz[baseline=-0.6ex]\draw[thick,orange,line width = 0.5mm] (0,0)--(0.5,0);}
\DeclareRobustCommand\fullc  {\tikz[baseline=-0.6ex]\draw[thick,orange,line width = 0.5mm] (0,0)--(0.5,0);}
\DeclareRobustCommand\fulld  {\tikz[baseline=-0.6ex]\draw[thick,orange,line width = 0.5mm] (0,0)--(0.5,0);}


\DeclareRobustCommand\full  {\tikz[baseline=-0.6ex]\draw[thick,orange,line width = 0.5mm] (0,0)--(0.5,0);}
\DeclareRobustCommand\dashed{\tikz[baseline=-0.6ex]\draw[thick,dashed, black,line width = 0.5mm] (0,0)--(0.54,0);}

\DeclareRobustCommand\fulla  {\tikz[baseline=-0.6ex]\draw[thick,red,line width = 0.5mm] (0,0)--(0.5,0);}

\DeclareRobustCommand\fullb  {\tikz[baseline=-0.6ex]\draw[thick,blue,line width = 0.5mm] (0,0)--(0.5,0);}

\begin{figure}
	\begin{center}
		\hspace{-1.cm}
		\includegraphics{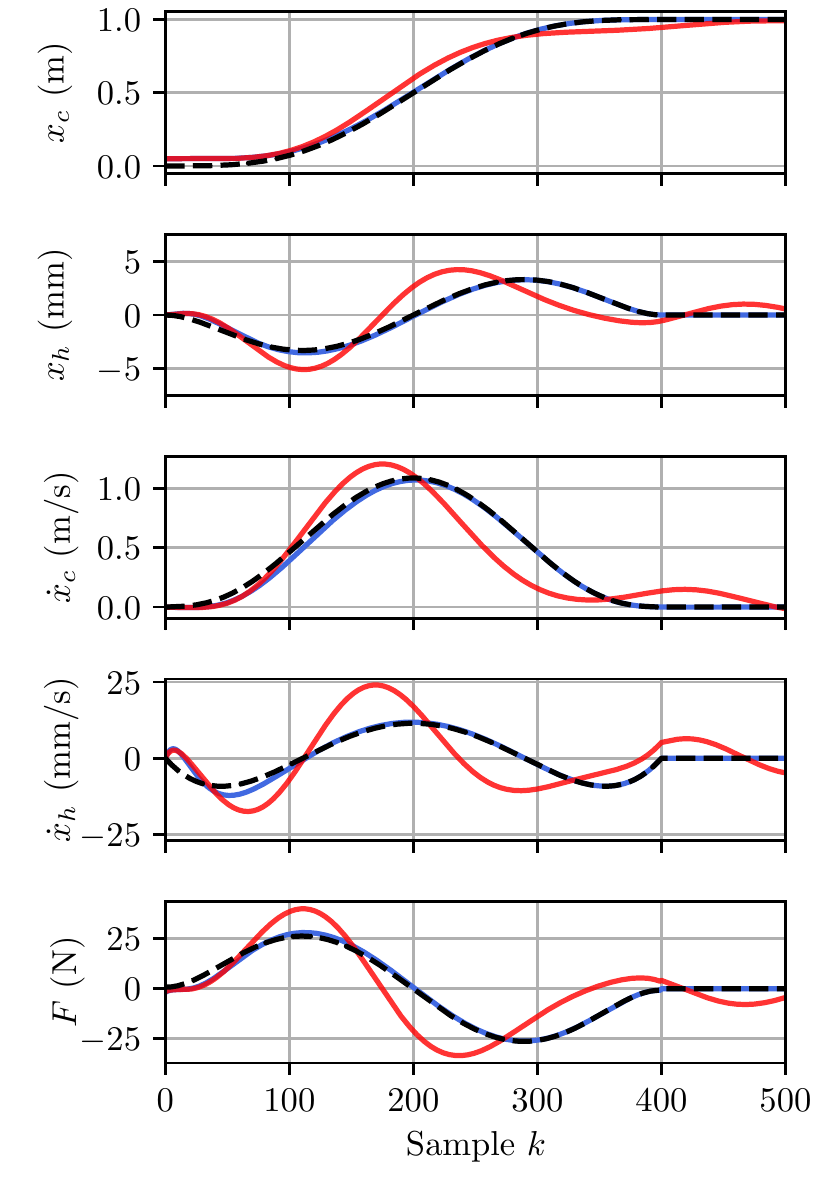}
		\caption{Comparison of state and input trajectories of the nominal controller applied to the crane $\Sigma_T \; (\fulla)$, the results (\fullb) after the controller has been learned to the new crane, and the desired trajectory (\dashed).}  
		\label{fig:smc_trajectory_04}
	\end{center}
\end{figure}


\section{Conclusion and Outlook}

In this paper, a data-driven control and transfer learning approach has been proposed and demonstrated in two examples.  
The idea  is to train a neural auto-encoder structure to represent the original system in a canonical Brunosvky structure and then to apply a standard controller design from linear control theory.
The design was validated in simulations on a pair of single-mast stacker cranes, where a nominal controller was learned on a nominal crane and then transferred to a crane identical in construction, but with modified parameters.
Based on the indirect data-driven feedback linearization and the resulting discrete-time Brunovsky form, trajectories were planned and stabilized for a change in operating points.
The transferred controller was able to improve its performance and reduce the need for manual adjustments after a relearning phase, where even occurring vibrations could be suppressed.
In terms of future work, the authors of this paper highlighted several directions that will be explored, such as extending the approach to a single-mast stacker crane with lifting unit leading to a MIMO system. 
In this context, the idea is to learn the representation of a multidimensional canonical structure. For this purpose, it will be necessary to deduce the lengths of the shift register chains on the basis of historical sensor data, in order to avoid the chain lengths as an additional hyperparameter.
Furthermore, the assumption that $\dim(z_k) = \dim(x_k)$ restricts this approach to input-state linearizable systems, which can be extended to the identification of flat systems by allowing $\dim(z_k) > \dim(x_k)$. Similarly, the authors also plan to investigate different and modified neural network structures to accelerate and improve the training process, such as incorporating a priori model information or bypassing the transformation constraints of the auto-encoder.


{\small
	\bibliographystyle{IEEEtranS}
	\bibliography{IEEEexample}

\begin{thebibliography}{10}
\providecommand{\url}[1]{#1}
\csname url@rmstyle\endcsname
\providecommand{\newblock}{\relax}
\providecommand{\bibinfo}[2]{#2}
\providecommand\BIBentrySTDinterwordspacing{\spaceskip=0pt\relax}
\providecommand\BIBentryALTinterwordstretchfactor{4}
\providecommand\BIBentryALTinterwordspacing{\spaceskip=\fontdimen2\font plus
\BIBentryALTinterwordstretchfactor\fontdimen3\font minus
  \fontdimen4\font\relax}
\providecommand\BIBforeignlanguage[2]{{%
\expandafter\ifx\csname l@#1\endcsname\relax
\typeout{** WARNING: IEEEtran.bst: No hyphenation pattern has been}%
\typeout{** loaded for the language `#1'. Using the pattern for}%
\typeout{** the default language instead.}%
\else
\language=\csname l@#1\endcsname
\fi
#2}}

\bibitem{Alsalti2021DataBasedSA}
M.~Alsalti, J.~Berberich, V.~G. Lopez, F.~Allg{\"o}wer, and M.~A. M{\"u}ller,
  ``Data-based system analysis and control of flat nonlinear systems,''
  \emph{2021 60th IEEE Conference on Decision and Control (CDC)}, pp.
  1484--1489, 2021.

\bibitem{Barrett2010TransferLF}
S.~Barrett, M.~E. Taylor, and P.~Stone, ``Transfer learning for reinforcement
  learning on a physical robot,'' in \emph{Adaptive Agents and Multi-Agent
  Systems}, 2010.

\bibitem{diwold2}
J.~Diwold, B.~Kolar, and M.~Schöberl, ``Discrete-time flatness-based control
  of a gantry crane,'' \emph{Control Engineering Practice}, vol. 119, 2 2022.

\bibitem{diwold1}
------, ``A trajectory-based approach to discrete-time flatness,'' \emph{IEEE
  Control System Letters}, vol.~6, pp. 289--294, 2022.

\bibitem{EckerMalzerWahrburg}
L.~Ecker, T.~Malzer, A.~Wahrburg, and M.~Schöberl, ``Observer design for a
  single mast stacker crane,'' \emph{at - Automatisierungstechnik}, vol.~69,
  no.~9, pp. 806--816, 2021.

\bibitem{Pasqualetti}
D.~Gadginmath, V.~Krishnan, and F.~Pasqualetti, ``Data-driven feedback
  linearization using the koopman generator,'' 2022.

\bibitem{Grizzle}
J.~W. Grizzle, ``Feedback linearization of discrete-time systems,'' in
  \emph{Analysis and Optimization of Systems}.\hskip 1em plus 0.5em minus
  0.4em\relax Springer Berlin Heidelberg, 1986, pp. 273--281.

\bibitem{isidori1995nonlinear}
A.~Isidori, \emph{Nonlinear Control Systems}, ser. Communications and Control
  Engineering.\hskip 1em plus 0.5em minus 0.4em\relax Springer London, 1995.

\bibitem{Krishnan2021}
V.~Krishnan and F.~Pasqualetti, ``On direct vs indirect data-driven predictive
  control,'' 12 2021, pp. 736--741.

\bibitem{BruntonAutoEncoder}
B.~Lusch, J.~Kutz, and S.~Brunton, ``Deep learning for universal linear
  embeddings of nonlinear dynamics,'' \emph{Nature Communications}, vol.~9, 11
  2018.

\bibitem{9328227}
B.~Maschler and M.~Weyrich, ``Deep transfer learning for industrial automation:
  A review and discussion of new techniques for data-driven machine learning,''
  \emph{IEEE Industrial Electronics Magazine}, vol.~15, no.~2, pp. 65--75,
  2021.

\bibitem{Nijmeijer_Nonlinear_dynamical_control_systems}
H.~Nijmeijer and A.~{Schaft, van der}, \emph{Nonlinear dynamical control
  systems}, corrected 2nd print.~ed.\hskip 1em plus 0.5em minus 0.4em\relax
  Springer, 1991.

\bibitem{7954034}
H.~Rams, M.~Schöberl, and K.~Schlacher, ``Optimal motion planning and
  energy-based control of a single mast stacker crane,'' \emph{IEEE
  Transactions on Control Systems Technology}, vol.~26, no.~4, pp. 1449--1457,
  2018.

\bibitem{Suykens1995ArtificialNN}
J.~A.~K. Suykens, J.~Vandewalle, and B.~D. Moor, ``Artificial neural networks
  for modelling and control of non-linear systems,'' 1995.

\bibitem{westenbroek2020feedback}
T.~Westenbroek, D.~Fridovich-Keil, E.~Mazumdar, S.~Arora, V.~Prabhu, S.~S.
  Sastry, and C.~J. Tomlin, ``Feedback linearization for uncertain systems via
  reinforcement learning,'' in \emph{2020 IEEE International Conference on
  Robotics and Automation (ICRA)}.\hskip 1em plus 0.5em minus 0.4em\relax IEEE,
  2020, pp. 1364--1371.

\bibitem{YESILDIREK19951659}
A.~Yeşildirek and F.~Lewis, ``Feedback linearization using neural networks,''
  \emph{Automatica}, vol.~31, no.~11, pp. 1659--1664, 1995.

\bibitem{Zahin}
S.~\'{z}Ahin, ``Learning feedback linearization using artificial neural
  networks,'' \emph{Neural Process. Lett.}, vol.~44, no.~3, p. 625–637, dec
  2016.

\end{thebibliography}
}

\newpage

\end{document}